\documentclass[11pt]{article}
\usepackage{amsfonts,latexsym,rawfonts,amsmath,amssymb,amsthm,graphicx}
\newcommand{\bthm}[2]{\vskip 8pt\noindent\bf #1\hskip 2pt\bf#2\it \hskip 8pt}
\newcommand{\ethm}{\vskip 8pt\rm}
\numberwithin{equation}{section}
\newtheorem{theorem}{Theorem}[section]
\newtheorem{lem}[theorem]{Lemma}
\newtheorem{thm}[theorem]{Theorem}

\newtheorem{cor}[theorem]{Corollary}
\newtheorem{defi}[theorem]{Definition}
\newtheorem{rem}[theorem]{Remark}
\def\s{\,\,\,\,}
\def\lan{\langle}
\def\ran{\rangle}

\def\mv{1.7ex}

\def\R{\mathbb{R}}

\def\K{\mathcal{K}}

\title{On conformal surfaces of annulus type}
\author{Yong Luo\footnote{The author is supported by the DFG Collaborative Research Center SFB/Transregio 71.}}
\date{}
\begin{document}
\maketitle
\begin{abstract}
Let $a>b>0$ and $f$ be a conformal map from $B_a\setminus B_b\subseteq R^2$ into
$\R^n$, with $|\nabla f|^2=2e^{2u}$. Then $(e_1, e_2)$ with $e_1=e^{-u}\frac{\partial f}{\partial r},$ and $e_2=r^{-1}e^{-u}\frac{\partial f}{\partial\theta}$ is a moving frame on $f(B_a\setminus B_b)$. It satisfies the following equation
$$d\star\langle de_1, e_2\rangle=0,$$
where $\star$ is the Hodge star operator on $R^2$ with respect to the standard metric.

We will study the Dirichret energy of this frame and give some applications.
\end{abstract}

\section{Introduction}

Let $\Omega$ be a smooth bounded domain in $R^2$ and $f$ be a $W^{2,2}$ map from $\Omega$ to $R^n$, and
 $(e_1,e_2)\in W^{1,2}(\Omega,\R^n\times \R^n)$ be a
positively oriented basis of $f$.
We define
$$\K(e_1,e_2)=
\frac{\partial e_1}{\partial x^1} \frac{\partial e_2}{\partial x^2}-
\frac{\partial e_1}{\partial x^2} \frac{\partial e_2}{\partial
x^1}=:\nabla e_1\nabla^\bot e_2.$$
It is easy to check that
 $\K(e_1,e_2)$ is invariant under the
group action $U(2)$, that is for any
$$e_1'=e_1\cos\theta+e_2\sin\theta,\s e_2'=
-e_1\sin\theta+e_2\cos\theta,$$
where $\theta\in W^{1,2}(\Omega, R)$, we have
\begin{eqnarray}
\K(e_1',e_2')=\K(e_1,e_2),
\end{eqnarray}
see the appendix A for a proof. Hence we can  write $\K(X_f):=\K(e_1,e_2)$, where $X_f$ is the Gauss map of the surface $f(\Omega)$ defined from $\Omega$ to the Grassmannian $G(2,n)$ . Moreover, we have
$$K_fe^{2v}=\K(X_f)=\nabla e_1\nabla^\bot
e_2,$$
where $K_f$ is the gauss curvature of the immersed surface $f(\Omega)$ and $|\nabla f|=2e^{2v}$ (see appendix A for the proof).

 Since $div\nabla^\bot
e_1=0$, and $rot\nabla e_2=0$, $\K(X_f)$ has compensation
compactness. Furthermore, Wente's type inequality can be applied here.

Recall Wente's type inequality, which states that if $a,b\in W^{1,2}(\Omega)$ and $u\in W^{1,2}_0(\Omega)$
solves the equation
$$-\Delta u=\nabla a\nabla^\bot b\s in\s \Omega,$$
then $u$ is continuous and we have
\begin{eqnarray}
\|u\|_{L^\infty(\Omega)}+\|\nabla u\|_{L^2(\Omega)}\leq C(\Omega)\|\nabla a\|_{L^2(\Omega)}
\|\nabla b\|_{L^2(\Omega)},
\end{eqnarray}
see (\cite{B},\cite{BC}).

It is easy to see that $C(\Omega)$ is invariant under translation and dilation. F.Bethuel and J.-M. Ghidaglia in \cite{BG1} showed that there exists constant $C_1$ which does not depend on $\Omega$ such that (1.2) holds true:
\begin{eqnarray*}
\|u\|_{L^\infty(\Omega)}+\|\nabla u\|_{L^2(\Omega)}\leq C_1\|\nabla a\|_{L^2(\Omega)}
\|\nabla b\|_{L^2(\Omega)}.
\end{eqnarray*}

We denote by $C_\infty(\Omega)$ the best constant involving the $L^\infty$ norm and by $C_2(\Omega)$ the best constant involving the $L^2$ norm. Then we have
\begin{eqnarray}
\|u\|_{L^\infty(\Omega)}&\leq& C_\infty(\Omega)\|\nabla a\|_{L^2(\Omega)}\|\nabla b\|_{L^2(\Omega)},
\\\|\nabla u\|_{L^2(\Omega)}&\leq& C_2(\Omega)\|\nabla a\|_{L^2(\Omega)}
\|\nabla b\|_{L^2(\Omega)}.
\end{eqnarray}

It is proved in \cite{B} that $C_\infty(\Omega)\geq\frac{1}{2\pi}$, and when $\Omega$ is simply connected, $C_\infty(\Omega)=\frac{1}{2\pi}$. For the general $\Omega$, it is proved by Topping \cite{Top} that $C_\infty(\Omega)=\frac{1}{2\pi}$. It is proved by Ge in \cite{G} that $C_2(\Omega)=\sqrt{\frac{3}{64\pi}}$.

Let $B\subseteq R^2$ be the unit disk centered at the origin, then Li, Luo and Tang proved the following theorem by using the inequality (1.3):

\bthm{Theorem}{ A\cite{L-L-T}} Let $\varphi\in W^{1,2}(B,G(2,n))$  with
$$\int_{B}|\mathcal{K}(\varphi)|d\sigma\leq\gamma<2\pi,$$
see appendix A for the information of $\mathcal{K}(\varphi)$.

Then there exists a map $(e_1,e_2) \in W^{1,2}(B,\R^n\times \R^n)$
such that for almost every $z\in B$, $(e_1(z),e_2(z))$ is a
positively oriented basis of $\varphi(z)$. Furthermore, we have
$$\|d(e_1,e_2)\|_{L^2(B)}\leq C(\gamma)\|\nabla \varphi\|_{L^2(B)}.$$
\ethm

 Note that $\K(\varphi)\leq\frac{1}{2}|\nabla \varphi|^2$, see (3.2). We have the following directly corollary:
 \bthm{Theorem}{ B\cite{L-L-T}} Let $\varphi\in W^{1,2}(B,G(2,n))$  with
$$\int_{B}|\nabla\varphi|^2d\sigma\leq\gamma<4\pi.$$

Then there exists a map $(e_1,e_2) \in W^{1,2}(B,\R^n\times \R^n)$
such that for almost every $z\in B$, $(e_1(z),e_2(z))$ is a
positively oriented basis of $\varphi(z)$. Furthermore, we have
$$\|d(e_1,e_2)\|_{L^2(B)}\leq C(\gamma)\|\nabla \varphi\|_{L^2(B)}.$$
\ethm
\begin{rem}
(1) The above theorem improved a theorem of H$\acute{e}$lein (\cite{H}, chapter 5) by changing the constant from $\frac{8\pi}{3}$ to $4\pi$ (the same result also is proved in \cite{K-L} by using a difficult result of \cite{M-V}). The difference between us is that in H$\acute{e}$lein's original proof he used the wente's inequality (1.4) of $L^2$ norm , whereas we use the wente's inequality (1.3) of $L^\infty$ norm.
\\(2) The constant $4\pi$ is shape for $n>3$ (see {\cite{K-L}}).
\end{rem}

Assume that $f: B\rightarrow R^n$ is a conformal map and $f\in C^\infty(\overline{B})$, and $X_f: f(B)\rightarrow G(2,n)$ is the gauss map. Let $\varphi=X_f\circ f\in  W^{1,2}(B,G(2,n))$, then we have (see (3.4) in the appendix A):
\begin{eqnarray*}
\int_B|K_f|d\mu_f&=&\int_B|\mathcal{K}(\varphi)|d\sigma,
\\ \int_B|\nabla\varphi|^2dx&=&\int_{f(B)}|\nabla_{g_f}X_f|^2d\mu_f=\int_B|A_f|^2d\mu_f,
\end{eqnarray*}
where $\nabla_{g_f}$ is the gradient with respect to $g_f$, $\mu_f$ is the area measure on $f(B)$, and $A_f$ is the second fundamental form of $f(B)$.

Hence for such a conformal immersion $f$ with the $L^1$ norm of the gauss curvature bounded by $2\pi$, there exists a moving frame on it whose Dirichlet energy is bounded by the $L^2$ norm of the second fundamental form. H$\acute{e}$lein (\cite{H}, chapter 5) used this moving frame to derive the weak compactness of immersed conformal surfaces from $B$ into $R^n$. Using his argument, we have
\bthm{Theorem}{C \cite{L-L-T}}
Let $f_k\in C^\infty(\bar{B},\R^n)$
be a sequence of conformal immersions with
$$\sup_k\int_B|K_{f_k}|d\mu_{f_k}\leq\gamma<2\pi,\s\sup_k\int_{B}|A_{f_k}|^2d\mu_{f_k}<\infty,$$
 where $d\mu_{f_k}$ is the volume form deduced from metric $g_{f_k}$.
Assume that $f_k$ converges to $f_0$  weakly in $W^{1,2}$. Then $f_0$ is
either a point or a conformal immersion.
\ethm

In this paper, we are interested in generalizing these above results to immersed conformal surfaces from $\Omega$ into $R^n$ when $\Omega$ is not simply connected. We will consider the easiest case, that is when $\Omega$ is an annuli. In the following we will let $a>b>0$, and $B_a\setminus B_b=\{x\in R^2: b<|x|<a\}$.
\begin{thm}
For every conformal map $f: B_a\setminus B_b\rightarrow R^n$
satisfying
$${\|\mathcal{K}(X_f)\|_{L^1(B_a\setminus B_b)}}\leq\gamma <2\pi,$$
there exists a map $b=(e_1,e_2)$ in $W^{1,2}(B_a\setminus B_b,\R^n\times R^n)$, such that for
almost every $z\in B_a\setminus B_b$, $(e_1(z),e_2(z))$ is a positively
 oriented basis of
$\varphi(z)$ and $\|d(e_1,e_2)\|_{L^2(B_a\setminus B_b)}$ is bounded.

Furthermore, if
\begin{eqnarray}
\frac{\beta}{1-\sqrt{\frac{\gamma}{2\pi}}}<1,
\end{eqnarray}
where
$$\int_{B_a\setminus B_b}|\frac{\partial e_2}{\partial \theta}d\theta|^2+|\frac{\partial e_1}{\partial \theta}d\theta|^2dx=\beta^2(\|\nabla e_1\|_{L^2(B_a\setminus B_b)}^2+\|\nabla
e_2\|_{L^2(B_a\setminus B_b)}^2).$$
then we have that
$$\|\nabla e_1\|_{L^2(B_a\setminus B_b)}^2+\|\nabla
e_2\|_{L^2(B_a\setminus B_b)}^2\leq C(\beta,\gamma,\frac{a}{b})\|A_f\|_{L^2(B_a\setminus B_b)}^2,$$
where $C(\beta,\gamma,\frac{a}{b})$ is a constant depending on $\beta$, $\gamma$ and $\frac{a}{b}$.
\end{thm}

As a direct corollary, we have
\begin{thm}
For every conformal map $f: B_a\setminus B_b\rightarrow R^n$
satisfying
$${\|A_f\|^2_{L^2(B_a\setminus B_b)}}\leq\gamma <4\pi,$$
where $A_f$ is the second fundamental form of $f$.

There exists a map $b=(e_1,e_2)$ in $W^{1,2}(B_a\setminus B_b,\R^n\times R^n)$, such that for
almost every $z\in B_a\setminus B_b$, $(e_1(z),e_2(z))$ is a positively
 oriented basis of
$f(z)$ and $\|d(e_1,e_2)\|_{L^2(B_a\setminus B_b)}$ is bounded.

Furthermore, if
\begin{eqnarray*}
\frac{\beta}{1-\sqrt{\frac{\gamma}{2\pi}}}<1,
\end{eqnarray*}
where
$$\int_{B_a\setminus B_b}|\frac{\partial e_2}{\partial \theta}d\theta|^2+|\frac{\partial e_1}{\partial \theta}d\theta|^2dx=\beta^2(\|\nabla e_1\|_{L^2(B_a\setminus B_b)}^2+\|\nabla
e_2\|_{L^2(B_a\setminus B_b)}^2).$$
then we have that
$$\|\nabla e_1\|_{L^2(B_a\setminus B_b)}^2+\|\nabla
e_2\|_{L^2(B_a\setminus B_b)}^2\leq C(\beta,\gamma,\frac{a}{b})\|A_f\|_{L^2(B_a\setminus B_b)}^2,$$
where $C(\beta,\gamma,\frac{a}{b})$ is a constant depending on $\beta$, $\gamma$ and $\frac{a}{b}$.
\end{thm}
\begin{rem}
We can see in the proof of the above two theorems that these estimates hold ture for any moving frame $(e_1, e_2)$ satisfying
$$d\star\langle e_1, de_2\rangle=0\s in\s B_a\setminus B_b,\s and\s\langle e_1, de_2\rangle(\frac{\partial}{\partial\nu})=0\s on\s \partial (B_a\setminus B_b),$$
where $\star$ is the Hodge star operator on $R^2$ and $\frac{\partial}{\partial\nu}$ is the outward normal vector on the boundaries. In the following we will define such a moving frame to be a coulomb frame.
\end{rem}
It is nature to ask the following question:
\\\textbf{Question 1:} On which kind of conformal parametric surfaces from $B_a\setminus B_b$ to $R^n$, there exists a moving frame $(e_1, e_2)$ on it and some $\beta\in (0, 1)$, such that
$$\int_{B_a\setminus B_b}|\frac{\partial e_2}{\partial \theta}d\theta|^2+|\frac{\partial e_1}{\partial \theta}d\theta|^2dx=\beta^2(\|\nabla e_1\|_{L^2(B_a\setminus B_b)}^2+\|\nabla
e_2\|_{L^2(B_a\setminus B_b)}^2)\textbf{?}$$

\begin{defi}
Let $\Omega$ be a domain in $R^2$ and $f: \Omega\rightarrow R^n$ be a conformal immersion, and $(e_1, e_2)$ is a moving frame on $f(\Omega)$, then we call  $(e_1, e_2)$ to be a coulomb frame of $f(\Omega)$ if
$$d\star\langle de_2, e_1\rangle=0\s in\s \Omega,\s \langle de_2, e_1\rangle(\frac{\partial}{\partial\nu})=0\s on\s \partial\Omega.$$
If we only have
$$d\star\langle de_2, e_1\rangle=0\s in\s \Omega,$$
then $(e_1, e_2)$ is called a semi-coulomb frame.
\end{defi}

We have
\begin{lem}
 Let $f: B_a\setminus B_b\rightarrow R^n$ be a conformal immersion, with $|\nabla f|^2=2e^{2u}$, and let $e_1=e^{-u}\frac{\partial f}{\partial r}$ and $e_2=r^{-1}e^{-u}\frac{\partial f}{\partial\theta}$, then $(e_1, e_2)$ is a semi-coulomb frame on $f(B_a\setminus B_b)$, and it is a coulomb frame if and only if u are constants on the boundaries. we call $(e_1, e_2)$ to be the canonical semi-coulomb frame.
\end{lem}
\begin{thm}
Let $f: B_a\setminus B_b\rightarrow R^n$ be a conformal minimal immersion, that is, $f(B_a\setminus B_b)$ is a conformal minimal surface in $R^n$, then the canonical semi-coulomb frame of $f(B_a\setminus B_b)$, $(e_1, e_2)$ satisfies the following condition
$$\int_{B_a\setminus B_b}|\frac{\partial e_2}{\partial \theta}d\theta|^2+|\frac{\partial e_1}{\partial \theta}d\theta|^2dx=\frac{1}{2}\int_{B_a\setminus B_b}|\nabla e_1|^2+|\nabla e_2|^2dx,$$
if and only if
$$\int_{B_a\setminus B_b}(\frac{1}{r}+\frac{\partial u}{\partial r})^2dx=\int_{B_a\setminus B_b}r^{-2}(\frac{\partial u}{\partial \theta})^2dx,$$
where $|\nabla f|^2=2e^{2u}.$
\end{thm}
As a direct corollary we have
\begin{cor}
Let $f$ be the same as the above theorem, and $u$ is radially symmetric, and let $(e_1, e_2)$ be the canonical semi-coulomb frame with
$$\int_{B_a\setminus B_b}|\frac{\partial e_2}{\partial \theta}d\theta|^2+|\frac{\partial e_1}{\partial \theta}d\theta|^2dx=\frac{1}{2}\int_{B_a\setminus B_b}|\nabla e_1|^2+|\nabla e_2|^2dx,$$
then we have
$$u(r)=-\log r+c,$$
where $c$ is some constant, and so $A=0$. Furthermore, $e_1$, $e_2$ are constant vectors.
\end{cor}
\begin{rem}
(1) Let $f$ be defined on $B_a\setminus B_b$ as $f(r,\theta)=(\pm e^c\log r,\theta,0,...,0)$, $c$ is a constant, then $f$ is a conformal immersion into $R^n$ with $u=-\log r+c$. It is easy to see that the canonical semi-coulomb frame of $f(B_a\setminus B_b)$ is ((1,0) (0,1)) and the second fundamental form of $f(B_a\setminus B_b)$ is zero.
\\ (2) Let $f$ be a conformal immersion from  $B_a\setminus B_b$ to $R^n$, with $|\nabla f|^2=2e^{2u}$, and $u$ is radially symmetric. Assume $(e_1,e_2)$ be the canonical semi-coulomb frame, then if $e_1$ and $e_2$ are constant vectors, then we have
$$\frac{\partial f}{\partial r}=e^u(a,b),\s \frac{\partial f}{\partial\theta}=re^u(c,d),$$
where $a, b, c, d$ are constants with $a^2+b^2=c^2+d^2=1, ac+bd=0$.
Thus $\frac{\partial^2f}{\partial r\partial\theta}=\frac{\partial^2f}{\partial\theta\partial r}$ implies that $u_r=-\frac{1}{r}$, and so $u(r)=-\log r+c$ for some constant $c$.
\end{rem}
As a corollary of the above results, we have the following theorem, which partially answers the above question:
\begin{thm}
Let $f: B_a\setminus B_b\rightarrow R^n$ be a conformal minimal immersion, with $|\nabla f|^2=2e^{2u}$, and
\begin{eqnarray}
\int_{B_a\setminus B_b}(\frac{1}{r}+\frac{\partial u}{\partial r})^2dx&=&\int_{B_a\setminus B_b}r^{-2}(\frac{\partial u}{\partial \theta})^2dx,
\\u&=&c_a\s on\s \partial B_a,\s u=c_b\s on\s \partial B_b,
\\\int_{B_a\setminus B_b}|K_f|d\mu_f&<&(3-2\sqrt{2})\pi,
\end{eqnarray}
where $c_a$ and $c_b$ are constants.
 Let $(e_1, e_2)$ is the canonical semi-coulomb frame of $f(B_a\setminus B_b)$, then we have
 \begin{eqnarray}
\int_{B_a\setminus B_b}|\nabla e_1|^2+|\nabla e_2|^2dx\leq C\int_{B_a\setminus B_b}|A_f|^2d\mu_f.
\end{eqnarray}
\end{thm}
\begin{rem}The property of coulomb frames on conformal surface $f(\Omega)$ has big difference between the case when $\Omega$ is simply connected and the case when $\Omega$ is not simply connected. Recall that if $f$ is a conformal immersion from the unit disk $B\subseteq R^2$ into $R^n$, then if the $L^2$ norm of the second fundamental form of $f(B)$ is below some constant, then the energy of a coulomb frame is controlled by this energy (\cite{H}\cite{L-L-T}). But for $f: B_a\setminus B_b\rightarrow R^n$ with $f(x,y)=(x,y,0,...,0)$, the second fundamental form of $f(B_a\setminus B_b)$ is zero, and $((\cos\theta,\sin\theta),(-\sin\theta,\cos\theta))$ is a coulomb frame on $f(B_a\setminus B_b)$ with nonzero energy.
\end{rem}
Now we will give an application of the above theorem.
\begin{thm}
Let $\{f_m\in C^\infty(\overline{B_a\setminus B_b})\}_{m\geq1}$ be a sequence of minimal conformal immersions from $B_a\setminus B_b$ into $R^n$ with $|\nabla f_m|^2=2e^{2u_m}$, and
\begin{eqnarray}
\int_{B_a\setminus B_b}(\frac{1}{r}+\frac{\partial u_m}{\partial r})^2dx&=&\int_{B_a\setminus B_b}r^{-2}(\frac{\partial u_m}{\partial \theta})^2dx,
\\u_m&=&c_{ma}\s on\s \partial B_a,\s u_m=c_{mb}\s on\s \partial B_b,
\\\sup_m\int_{B_a\setminus B_b}|K_{f_m}|d\mu_m&<&(3-2\sqrt{2})\pi,\s \sup_m\int_{B_a\setminus B_b}|A_{f_m}|^2d\mu_m<\infty,
\end{eqnarray}
where $c_{ma}$ and $c_{mb}$ are constants, with $\sup_m\{|c_{ma}|+|c_{mb}|\}<\infty$.
\\Assume that $f_m$ converges weakly to $f_0$ in $W^{1,2}$, then $f_0$ is a minimal conformal immersion, with bounded conformal factor. Furthermore the metric induced by $f_0$ is continuous.
\end{thm}
\begin{rem}
The minimal property will be kept under the weak convergence by the definition, if the limit immersion is conformal. Hence the difficult and non-obvious part is to prove that the limit immersion is conformal, and with bounded conformal factor.
\end{rem}
This paper is organized as follows: In section 2 results are proved. In appendix A we will give some basic computations which have been used in this paper and in appendix B we will give an alternative proof about that $f_0$ is conformal in theorem 1.12 by using a strong converge theorem of p-harmonic maps due to Hardt, Lin and Mou (\cite{HLM}).

\textbf{Notations:} $\partial_r=\frac{\partial}{\partial r}$, $\partial_\theta=\frac{\partial}{\partial_\theta}$, $f_r=\frac{\partial f}{\partial r}$, $f_\theta=\frac{\partial f}{\partial_\theta}$, $\partial_{rr}=...etc.$

\textbf{Acknowledgment} The author would like to thank his advisor, professor Guofang Wang, and professor Ernst Kuwert for discussions.

\section{Proof of the results}

\textbf{Proof of theorem 1.2:}

Let $(e_1,e_2)\in
W^{1,2}(B_a\setminus B_b,\R^n\times \R^n)$ be a positively oriented basis of
$f$ and $X_f$ is the gauss map of $f$. Let
$$\begin{array}{ll}
    \digamma=\{(X_f,e_1,e_2)\in G(2,n)\times\R^n\times\R^n|
&(e_1,e_2)\\
     &\hbox{
is a positively oriented orthonormal  basis of }X_f\}.
\end{array}$$
Then $\digamma$ is a fibre bundle over $G(2,n)$ with fibre $S^1$.
Since $f$ is conformal, there exists a section $(\widetilde{e_1},
\widetilde{e_2})$ of $X_f^*\digamma$.

We consider for each $\theta\in W^{1,2}(B_a\setminus B_b,\R)$ the frame $(e_1,e_2)$
obtained by
$$ (e_1, e_2)=(\widetilde{e_1},
\widetilde{e_2}) \left(\begin{array}{cc}
\cos\theta & -\sin\theta\\
\sin\theta & \cos\theta
\end{array}\right).
$$
Actually, we will minimize over  $\theta\in W^{1,2}(B_a\setminus B_b,\R) $
the   functional
\begin{eqnarray}
 F(\theta)&=&\frac{1}{2}\int_{B_a\setminus B_b}(|\nabla e_1|^2+|\nabla
 e_2|^2)d\sigma\\&=&\int_{B_a\setminus B_b}|\omega_2^1|^2d\sigma,
\end{eqnarray}
where $w_2^1=\lan de_1,e_2\ran$. By the arguments in \cite{H}, the
minimum of F is attained, and  the minimizer $(e_1,e_2)$ satisfies
$$\left\{\begin{array}{rl}
d(\star\omega_2^1)= 0& in\s B_a\setminus B_b,\\[\mv]
{\omega_2^1(\frac{\partial}{\partial{\nu}}) }=0&on\s {\partial(B_a\setminus B_b)},
\end{array}\right.$$
where $\star$ is the Hodge star operator and $\frac{\partial}{\partial\nu}$ is the outward normal vector on the boundary.
Then there exists some $v\in
W^{1,2}(B_a\setminus B_b,\R)$ such that
\begin{equation}\label{dv}
  dv=\star\omega_2^1-\alpha d\theta\s in\s  B_a\setminus B_b,
\end{equation}
where $\alpha=\frac{1}{2\pi}\int_0^{2\pi}\star\omega_2^1$ is a constant.
It is easy to check that $\frac{\partial v}{\partial \theta}=
dv(\frac{\partial}{\partial{\theta}})=-\alpha$ on $\partial (B_a\setminus B_b)$, and hence we have
 $v|_{\partial B_a}=c_a-\alpha \theta,$ and $v|_{\partial B_b}=c_b-\alpha \theta,$ where $c_a=v(a,0)$ and $c_b=v(b,0).$
A direct calculation yields
\begin{equation}\label{eq.v}
{-\triangle}v=\K(e_1,e_2).
\end{equation}
Decompose $v$ to be $v=v_1+v_2$ where
\begin{equation}\label{eq.v}
{-\triangle}v_1=\K(e_1,e_2), v_1=0\s on\s \partial(B_a\setminus B_b).
\end{equation}
and
\begin{equation}\label{eq.v}
{-\triangle}v_2=0, v_2=c_1-\alpha\theta\s on\s \partial B_a, v_2=c_2-\alpha\theta\s on\s \partial B_b.
\end{equation}
To estimate $v_1$, we set $v_1^k$
to be the solution of
$$-\Delta v_1^k=\frac{\partial e_1^k}
{\partial x^1}\frac{\partial e_2^k} {\partial x^2}-\frac{\partial
e_1^k} {\partial x^2}\frac{\partial e_2^k} {\partial x^1},\s
v_1^k|_{\partial (B_a\setminus B_b)}=0,$$
where $e_i=(e_i^1,e_i^2,\cdots, e_i^n)$.
 Applying Wente's inequality, we have
$$\|v_1^k\|_{L^\infty(B_a\setminus B_b)}\leq \frac{1}{2\pi}
\|\nabla e_1^k\|_{L^2(B_a\setminus B_b)}\|\nabla e_2^k\|_{L^2(B_a\setminus B_b)},$$
which obviously implies that
$$\|v_1\|_{L^\infty(B_a\setminus B_b)}\leq\frac{1}{4\pi} (\|\nabla e_1\|_{L^2(B_a\setminus B_b)}^2+\|\nabla
e_2\|_{L^2(B_a\setminus B_b)}^2).$$
A simple calculation by integration by parts implies that
\begin{eqnarray*}
\int_{B_a\setminus B_b}|\nabla
v_1|^2d\sigma&=&\int_{B_a\setminus B_b}v_1\K(e_1,e_2)d\sigma\\
&\leq&\|v_1\|_{L^\infty(B_a\setminus B_b)}\int_{B_a\setminus B_b}|\K(e_1,e_2)|d\sigma\\
&\leq&\frac{\gamma }{4\pi}  (\|\nabla e_1\|_{L^2(B_a\setminus B_b)}^2+\|\nabla
e_2\|_{L^2(B_a\setminus B_b)}^2).
\end{eqnarray*}
For $v_2$, we have
$$v_2=\frac{c_a-c_b}{\log\frac{a}{b}}\log|x|+\frac{c_b\log a-c_a\log b}{\log\frac{a}{b}}-\alpha\theta.$$
Note that by the calculation in appendix A section (see (3.1)) we have
$$\|\nabla e_1\|_{L^2(B_a\setminus B_b)}^2+\|\nabla
e_2\|_{L^2(B_a\setminus B_b)}^2=2\|dv+\alpha d\theta\|_{L^2(B_a\setminus B_b)}^2+\|\nabla X_f\|_{L^2(B_a\setminus B_b)}^2.$$
We have
\begin{eqnarray*}
\|dv+\alpha d\theta\|_{L^2(B_a\setminus B_b)}&=&\|dv_1+dv_2+\alpha d\theta\|_{L^2(B_a\setminus B_b)}
\\&\leq&\|dv_1\|_{L^2(B_a\setminus B_b)}+\|dv_2+\alpha d\theta\|_{L^2(B_a\setminus B_b)}
\\&\leq&(\frac{\gamma}{4\pi})^\frac{1}{2}(\|\nabla e_1\|_{L^2(B_a\setminus B_b)}^2+\|\nabla
e_2\|_{L^2(B_a\setminus B_b)}^2)^\frac{1}{2}+(2\pi)^\frac{1}{2}\frac{|c_a-c_b|}{(\log\frac{a}{b})^\frac{1}{2}}.
\end{eqnarray*}
Thus we obtain
\begin{eqnarray}
(1-(\frac{\gamma}{2\pi})^\frac{1}{2})(\|\nabla e_1\|_{L^2(B_a\setminus B_b)}^2+\|\nabla
e_2\|_{L^2(B_a\setminus B_b)}^2)^\frac{1}{2}\leq(4\pi)^\frac{1}{2}\frac{|c_a-c_b|}{(\log\frac{a}{b})^\frac{1}{2}}+\|\nabla X_f\|_{L^2(B_a\setminus B_b)}.
\end{eqnarray}
Without loss of generality, we assume that $\int_b^a|\frac{\partial e_2}{\partial \theta}|^2(r,0)+|\frac{\partial e_1}{\partial \theta}|^2(r,0)dr\leq\int_b^a|\frac{\partial e_2}{\partial \theta}|^2(r,\theta)+|\frac{\partial e_1}{\partial \theta}|^2(r,\theta)dr,$ for $0\leq\theta<2\pi$.
To estimate the number $|c_a-c_b|$, we note that
\begin{eqnarray*}
|c_a-c_b|=|\int_b^a\frac{\partial v}{\partial r}(r,0)dr|&\leq&\int_b^a|\langle e_1\frac{\partial e_2}{\partial \theta}\rangle|(r,0)r^{-1} dr
\\&\leq&(\frac{\log\frac{a}{b}}{2})^\frac{1}{2}(\int_b^a(|\langle e_1\frac{\partial e_2}{\partial \theta}\rangle|^2+|\langle e_2\frac{\partial e_1}{\partial \theta}\rangle|^2)r^{-1} dr)^\frac{1}{2}
\\&\leq&(\frac{\log\frac{a}{b}}{2})^\frac{1}{2}(\int_b^a|\frac{\partial e_2}{\partial \theta}d\theta|^2+|\frac{\partial e_1}{\partial \theta}d\theta|^2rdr)^\frac{1}{2}
\\&\leq&(\frac{\log\frac{a}{b}}{4\pi})^\frac{1}{2}(\int_{B_a\setminus B_b}|\frac{\partial e_2}{\partial \theta}d\theta|^2+|\frac{\partial e_1}{\partial \theta}d\theta|^2dx)^\frac{1}{2}
\\&\leq&\beta(\frac{\log\frac{a}{b}}{4\pi})^\frac{1}{2}(\|\nabla e_1\|_{L^2(B_a\setminus B_b)}^2+\|\nabla
e_2\|_{L^2(B_a\setminus B_b)}^2)^\frac{1}{2}
\\&\leq&\frac{\beta}{1-(\frac{\gamma}{2\pi})^\frac{1}{2}}|c_a-c_b|+C(\gamma,\frac{a}{b})\|\nabla X_f\|_{L^2(B_a\setminus B_b)},
\end{eqnarray*}
where we have used the fact that $\star d\theta=r^{-1}dr$ and that
 $$\int_{B_a\setminus B_b}|\frac{\partial e_2}{\partial \theta}d\theta|^2+|\frac{\partial e_1}{\partial \theta}d\theta|^2dx=\beta^2(\|\nabla e_1\|_{L^2(B_a\setminus B_b)}^2+\|\nabla
e_2\|_{L^2(B_a\setminus B_b)}^2).$$

Hence if
$$\frac{\beta}{1-\sqrt{\frac{\gamma}{2\pi}}}<1,$$
 then $|c_a-c_b|$ is controlled by $C(\beta,\gamma,\frac{a}{b})\|\nabla X_f\|_{L^2(B_a\setminus B_b)},$ and hence the energy of the coulomb frame is controlled by $C(\beta,\gamma,\frac{a}{b})\|\nabla X_f\|_{L^2(B_a\setminus B_b)}.$  Noting that $\|\nabla X_f\|_{L^2(B_a\setminus B_b)}=\|A_f\|_{L^2(B_a\setminus B_b)},$ we complete the proof of theorem 1.2 .

 \textbf{Proof of lemma 1.6:}
 Note that $e_1=e^{-u}f_r$, and $e_2=r^{-1}e^{-u}f_\theta$, then we have
 \begin{eqnarray*}
 \langle de_1, e_2\rangle&=&\langle\frac{f_{rr-f_ru_r}}{e^u},\frac{f_\theta}{re^u}\rangle dr+\langle\frac{f_{r\theta-f_ru_\theta}}{e^u},\frac{f_\theta}{re^u}\rangle d\theta
 \\&=&\langle\frac{f_{rr}}{e^u},\frac{f_\theta}{re^u}\rangle dr+\langle\frac{f_{r\theta}}{e^u},\frac{f_\theta}{re^u}\rangle d\theta
 \\&=&\frac{-u_\theta}{r}dr+(1+ru_r)d\theta,
 \end{eqnarray*}
 hence
 \begin{eqnarray*}
 \star\langle de_1, e_2\rangle&=&\frac{-u_\theta}{r}dr+(1+ru_r)d\theta
 \\&=&\frac{-u_\theta}{r}(-rd\theta)+(1+ru_r)r^{-1}dr
 \\&=&u_\theta d\theta+(r^{-1}+u_r)dr,
 \end{eqnarray*}
 finally we obtain
 \begin{eqnarray*}
 d \star\langle de_1, e_2\rangle=u_{\theta r}dr\wedge d\theta+u_{r\theta}d\theta\wedge dr=0,
\end{eqnarray*}
 which implies that $(e_1, e_2)$ is a semi-coulomb frame on $f(B_a\setminus B_b)$. In addition, $\star\langle de_2, e_1\rangle(\frac{\partial}{\partial n})=0$ if and only if $\star\langle de_1, e_2\rangle(\frac{\partial}{\partial n})=0$ if and only if  $\langle de_1, e_2\rangle(\frac{\partial}{\partial\theta})=0$. Note that $\langle de_1, e_2\rangle(\frac{\partial}{\partial\theta})=u_\theta$, so $\star\langle de_1, e_2\rangle(\frac{\partial}{\partial n})=0$ if and only if $u$ are constants on the boundaries.

 \textbf{Proof of theorem 1.7:}
 Recall that $e_1=e^{-u}\frac{\partial f}{\partial r}$ and $e_2=r^{-1}e^{-u}\frac{\partial f}{\partial\theta}$, then we have
 \begin{eqnarray*}
 \frac{\partial e_1}{\partial \theta}&=&\frac{f_{r\theta}-f_ru_\theta}{e^u},
 \\ \frac{\partial e_2}{\partial\theta}&=&\frac{f_{\theta\theta}-f_\theta u_\theta}{re^u},
 \end{eqnarray*}
 hence we can obtain
 \begin{eqnarray*}
 |\frac{\partial e_1}{\partial \theta}|^2=&&e^{-2u}[f_{r\theta}^2-2\langle f_{r\theta}, f_r\rangle u_\theta+f_r^2u_\theta^2]
 \\&=&e^{-2u}[f_{r\theta}^2-2e^{2u}u_\theta^2+e^{2u}u_\theta^2]
 \\&=&e^{-2u}f_{r\theta}^2-u_\theta^2.
 \end{eqnarray*}
 and
 \begin{eqnarray*}
 |\frac{\partial e_2}{\partial \theta}|^2=&&r^{-2}e^{-2u}[f_{\theta\theta}^2-2\langle f_{\theta\theta}, f_\theta\rangle u_\theta+f_\theta^2u_\theta^2]
 \\&=&r^{-2}e^{-2u}[f_{\theta\theta}^2-\frac{\partial r^2e^{2u}}{\partial\theta}u_\theta+f_\theta^2u_\theta^2]
 \\&=&r^{-2}e^{-2u}[f_{\theta\theta}^2-2r^2e^{2u}u_\theta^2+r^2e^{2u}u_\theta^2]
 \\&=&r^{-2}e^{-2u}f_{\theta\theta}^2-u_\theta^2.
 \end{eqnarray*}
 Similarly, we can obtain
 \begin{eqnarray*}
  \frac{\partial e_1}{\partial r}&=&\frac{f_{rr}-f_ru_r}{e^u},
 \\ \frac{\partial e_2}{\partial r}&=&\frac{f_{r\theta}-f_\theta(\frac{1}{r}+ u_r)}{re^u},
 \end{eqnarray*}
 and
 \begin{eqnarray*}
 |\frac{\partial e_1}{\partial r}|^2&=&e^{-2u}f_{rr}^2-u_r^2,
 \\ |\frac{\partial e_2}{\partial r}|^2&=&r^{-2}e^{-2u}f_{r\theta}^2-(\frac{1}{r}+u_r)^2.
 \end{eqnarray*}
 Summarized the above computations and note that $|d\theta|^2=r^{-2},$ we have
 \begin{eqnarray}
 |\frac{\partial e_2}{\partial \theta}d\theta|^2+|\frac{\partial e_1}{\partial \theta}d\theta|^2&=&r^{-2}e^{-2u}f_{r\theta}^2+r^{-4}e^{-2u}f_{\theta\theta}^2-2r^{-2}u_\theta^2,
 \\|\frac{\partial e_2}{\partial r}dr|^2+|\frac{\partial e_1}{\partial r}dr|^2&=&e^{-2u}f_{rr}^2+r^{-2}e^{-2u}f_{r\theta}^2-u_r^2-(\frac{1}{r}+u_r)^2.
 \end{eqnarray}
 On the other hand, by the definition of the second fundamental form, we have
 \begin{eqnarray*}
 A_{rr}&=&f_{rr}-e^{-2u}\langle f_{rr}, f_r\rangle f_r-r^{-2}e^{-2u}\langle f_{rr}, f_\theta\rangle f_\theta,
 \\A_{\theta\theta}&=&f_{\theta\theta}-e^{-2u}\langle f_{\theta\theta}, f_r\rangle f_r-r^{-2}e^{-2u}\langle f_{\theta\theta}, f_\theta\rangle f_\theta,
 \end{eqnarray*}
 therefor we have
 \begin{eqnarray*}
 A_{rr}^2&=&f_{rr}^2-2\langle f_{rr}, f_r\rangle^2e^{-2u}-2\langle f_{rr}, f_\theta\rangle^2r^{-2}e^{-2u}+\langle f_{rr}, f_r\rangle^2e^{-2u}+\langle f_{rr}, f_\theta\rangle^2r^{-2}e^{-2u}
 \\&=&f_{rr}^2-\langle f_{rr}, f_r\rangle^2e^{-2u}-\langle f_{rr}, f_\theta\rangle^2r^{-2}e^{-2u}
 \\&=&f_{rr}^2-u_r^2e^{2u}-u_\theta^2r^{-2}e^{2u},
 \end{eqnarray*}
 similar computations implies that
 $$A_{\theta\theta}^2=f_{\theta\theta}^2-r^4e^{2u}(\frac{1}{r}+u_r)^2-r^2e^{2u}u_\theta^2.$$
 Note that $f$ is minimal, so
 $$Trace(A)=g^{rr}A_{rr}+2g^{r\theta}A_{r\theta}+g^{\theta\theta}A_{\theta\theta}=g^{rr}A_{rr}+g^{\theta\theta}A_{\theta\theta}=0,$$
 where $g^{rr}=e^{-2u}, g^{\theta\theta}=r^{-2}e^{-2u}$
 hence
 $$A_{rr}=-r^{-2}A_{\theta\theta},$$
 which implies that
 $$A_{rr}^2=r^{-4}A_{\theta\theta}^2.$$
 Thus we obtain
 \begin{eqnarray}
 f_{rr}^2-u_r^2e^{2u}=r^{-4}f_{\theta\theta}^2-e^{2u}(\frac{1}{r}+u_r)^2.
 \end{eqnarray}
 Note that
$$\int_{B_a\setminus B_b}|\frac{\partial e_2}{\partial \theta}d\theta|^2+|\frac{\partial e_1}{\partial \theta}d\theta|^2dx=\frac{1}{2}\int_{B_a\setminus B_b}|\nabla e_1|^2+|\nabla e_2|^2dx,$$
 if and only if
 $$\int_{B_a\setminus B_b}|\frac{\partial e_2}{\partial \theta}d\theta|^2+|\frac{\partial e_1}{\partial \theta}d\theta|^2dx=\int_{B_a\setminus B_b}|\frac{\partial e_2}{\partial r}dr|^2+|\frac{\partial e_1}{\partial r}dr|^2dx,$$
 thus we can get by combining (2.8)-(2.10) that
 $$\int_{B_a\setminus B_b}|\frac{\partial e_2}{\partial \theta}d\theta|^2+|\frac{\partial e_1}{\partial \theta}d\theta|^2dx=\frac{1}{2}\int_{B_a\setminus B_b}|\nabla e_1|^2+|\nabla e_2|^2dx,$$
 if and only if
  $$\int_{B_a\setminus B_b}(\frac{1}{r}+\frac{\partial u}{\partial r})^2dx=\int_{B_a\setminus B_b}r^{-2}(\frac{\partial u}{\partial \theta})^2dx.$$

\textbf{Proof of corollary 1.8:}
From theorem 1.7 we know that
$$\int_{B_a\setminus B_b}|\frac{\partial e_2}{\partial \theta}d\theta|^2+|\frac{\partial e_1}{\partial \theta}d\theta|^2dx=\frac{1}{2}\int_{B_a\setminus B_b}|\nabla e_1|^2+|\nabla e_2|^2dx,$$
implies
  $$\int_{B_a\setminus B_b}(\frac{1}{r}+\frac{\partial u}{\partial r})^2dx=\int_{B_a\setminus B_b}r^{-2}(\frac{\partial u}{\partial \theta})^2dx,$$
hence if $u$ is radially symmetric we must have that
$$\int_{B_a\setminus B_b}(\frac{1}{r}+\frac{\partial u}{\partial r})^2dx=0,$$
which implies that
$$\frac{1}{r}+\frac{\partial u}{\partial r}=0,$$
and so there is some constant $c$ such that
$$u(r)=-\log r+c.$$
Hence we have $\triangle u=0$, which implies that the gauss curvature $K=0$, and so $A=0$. By theorem 1.10 (note that $u$ are constants on boundaries and hence $(e_1, e_2)$ is a coulomb frame), we know that $(e_1, e_2)$ has zero energy and so $e_1$ and $e_2$ are constant vectors.

\textbf{Proof of theorem 1.10:} We know that $(e_1, e_2)$ is a coulomb frame and so when (1.6) holds we have that
$$\int_{B_a\setminus B_b}|\frac{\partial e_2}{\partial \theta}d\theta|^2+|\frac{\partial e_1}{\partial \theta}d\theta|^2dx=\frac{1}{2}\int_{B_a\setminus B_b}|\nabla e_1|^2+|\nabla e_2|^2dx,$$
by theorem 1.7. Then the constant $\beta$ in theorem 1.2 is $\frac{\sqrt{2}}{2}$ and so when
$$\int_{B_a\setminus B_b}|K_f|du_f<(3-2\sqrt{2})\pi,$$
we have (1.5) holds, and then we get the desired inequality (1.9) from theorem 1.2.

\textbf{Proof of theorem 1.12:}
Let $(e_{m1}, e_{m2})$ be the canonical semi-coulomb frame on $f_m(B_a\setminus B_b)$, then by theorem 1.10 we have the following inequality
\begin{eqnarray*}
\int_{B_a\setminus B_b}|\nabla e_{m1}|^2+|\nabla e_{m2}|^2dx\leq C\int_{B_a\setminus B_b}|A_m|^2d\mu_{f_m},
\end{eqnarray*}
where $C$ only depends on $\frac{a}{b}$.

Note that we have
$$-\Delta u_m=K_me^{2u_m}=\nabla e_{m1}\nabla^\bot e_{m2}\s in\s B_a\setminus B_b,$$
where $K_m$ is the gauss curvature and $\nabla=(\frac{\partial}{\partial x}, \frac{\partial}{\partial y})$ and $\nabla^\bot=(-\frac{\partial}{\partial y}, \frac{\partial}{\partial x})$.
Let $v_m$ solves the following
\begin{eqnarray*}
-\Delta v_m&=&\nabla e_{m1}\nabla^\bot e_{m2}\s in\s B_a\setminus B_b,
\\v_m&=&0\s on\s \partial(B_a\setminus B_b).
\end{eqnarray*}
Let $e_{mi}=(e_{mi}^1,...,e_{mi}^n),\s i=1,2$, and $v_m=v_m^1+...+v_m^n$, such that for each $1\leq k\leq n$,
\begin{eqnarray*}
-\Delta v_m^k&=&\nabla e_{m1}^k\nabla^\bot e_{m2}^k\s in\s B_a\setminus B_b,
\\v_m^k&=&0\s on\s \partial(B_a\setminus B_b),
\end{eqnarray*}
then by wente's inequality we obtain
$$\|v_m^k\|_{L^\infty(B_a\setminus B_b)}\leq \frac{1}{2\pi}\|\nabla e_{m1}^k\|_{L^2(B_a\setminus B_b)}\|\nabla e_{m2}^k\|_{L^2(B_a\setminus B_b)},$$
hence
\begin{eqnarray*}
\|v_m\|_{L^\infty(B_a\setminus B_b)}&\leq&\sum_k\|v_m^k\|_{L^\infty(B_a\setminus B_b)}
\\&\leq&\sum_k\frac{1}{2\pi}\|\nabla e_{m1}^k\|_{L^2(B_a\setminus B_b)}\|\nabla e_{m2}^k\|_{L^2(B_a\setminus B_b)}
\\&\leq&\frac{1}{2\pi}\|\nabla e_{m1}\|_{L^2(B_a\setminus B_b)}\|\nabla e_{m2}\|_{L^2(B_a\setminus B_b)},
\end{eqnarray*}
in the last inequality we have used Holder's inequality.

By using the equation satisfied by $v_m$ and by integration by parts we have
\begin{eqnarray*}
\int_{B_a\setminus B_b} -v_m\Delta v_m&=&\int v_m\nabla e_{m1}\nabla^\bot e_{m2}dx
\\&\leq&\|v_m\|_{L^\infty(B_a\setminus B_b)}\|\nabla e_{m1}\|_{L^2(B_a\setminus B_b)}\|\nabla e_{m2}\|_{L^2(B_a\setminus B_b)}
\\&\leq&\frac{1}{2\pi}\|\nabla e_{m1}\|_{L^2(B_a\setminus B_b)}^2\|\nabla e_{m2}\|_{L^2(B_a\setminus B_b)}^2.
\end{eqnarray*}
That is
\begin{eqnarray*}
\|v_m\|_{L^2(B_a\setminus B_b)}^2\leq \frac{1}{2\pi}\|\nabla e_{m1}\|_{L^2(B_a\setminus B_b)}^2\|\nabla e_{m2}\|_{L^2(B_a\setminus B_b)}^2.
\end{eqnarray*}
On the other hand,
\begin{eqnarray*}
\Delta(u_m-v_m)&=&0\s in\s B_a\setminus B_b,
\\u_m-v_m&=&c_{ma}\s on\s \partial B_a,
\\ u_m-v_m&=&c_{mb}\s on\s \partial B_b,
 \end{eqnarray*}
 thus we have
 $$u_m-v_m=\frac{c_{ma}-c_{mb}}{\log\frac{a}{b}}\log|x|+\frac{c_{mb}\log a-c_{ma}\log b}{\log\frac{a}{b}},$$
 which implies that
 \begin{eqnarray}
\|u_m\|_{L^2(B_a\setminus B_b)} +\|u_m\|_{L^\infty(B_a\setminus B_b)}\leq C<\infty,
 \end{eqnarray}
 for some constant $C$ independent of $m$.

 Then by using an argument given by (\cite{H}, chapter 5), we can get that $f_0$ is a conformal immersion with bounded conformal factor as the following: Because $f_m$ is conformal, there exists $0\leq\theta_m\in C^\infty<2\pi$ such that
 \begin{eqnarray}
 df_m=e^{u_m}((\cos\theta_me_{m1}+\sin\theta_me_{m2})dx+(-\sin\theta_me_{m1}+\cos\theta_me_{m2})dy).
 \end{eqnarray}
 In particular, projecting the equation $d^2f_m=0$ along $e_{m1}$ and $e_{m2}$ we obtain
 \begin{eqnarray}
 \frac{\partial\theta_m}{\partial x}+\frac{\partial u_m}{\partial y}&=&\omega_{m2}^1(\frac{\partial}{\partial x}),
 \\\frac{\partial\theta_m}{\partial y}-\frac{\partial u_m}{\partial x}&=&\omega_{m2}^1(\frac{\partial}{\partial y}),
 \end{eqnarray}
 where $\omega_{m2}^1=\langle de_{m2}, e_{m1}\rangle$.

 Note that (2.13)-(2.14) implies that $\theta_m$ is bounded in $W^{1,2}$, hence we have that (we do not distinguish a sequence and its subsequences)
 $$(b_m,\theta_m,u_m)\rightharpoonup(b,\theta,u)\s weakly\s in\s W^{1,2},$$
 and so
 $$(b_m,\theta_m,u_m)\to (b,\theta,u)\s in\s L^2,$$
 therefor we have
  $$(b_m,\theta_m,u_m)\to (b,\theta,u)\s a.e.\s in\s B_a\setminus B_b,$$
  where $b_m=(e_{m1}, e_{m2}),$ and $b=(e_{1}, e_{2}).$

 By passing to the limit in (2.12) we get
   \begin{eqnarray}
 df_0=e^u((\cos\theta e_1+\sin\theta e_2)dx+(-\sin\theta e_1+\cos\theta e_2)dy),
 \end{eqnarray}
 which implies that $f_0$ is conformal, with bounded conformal factor $e^u$.

 Because u satisfies the following wente's type equation
 $$-\Delta u=\nabla e_1\nabla^\bot e_2\s in\s B_a\setminus B_b,$$
hence u is continuous.

  Note that
 $$\Delta f_m=0\s in\s B_a\setminus B_b,$$
 and
 $$f_m\rightarrow f_0\s weakly\s in\s W^{1,2}(B_a\setminus B_b),$$
 therefor we have that
 \begin{eqnarray}
 \Delta f_0=0.
 \end{eqnarray}
 On the other hand, because $f_0$ is a conformal immersion with $|\nabla f_0|^2=2e^{2u}$, we have that
  \begin{eqnarray}
 \Delta f_0=e^{2u}H_{f_0},
 \end{eqnarray}
 where $H_{f_0}$ is the mean curvature vector of $f_0$.

 By comparing (2.16) with (2.17), we get that  $H_{f_0}=0$, and so $f_0$ is a minimal immersion.

\section{Appendix}

\subsection{A}
In this appendix, we review briefly  some basic facts of
Grassmannian. The concept in this appendix can be found in any
textbook on the theory of Grassmannian.

Let
$$\Lambda^{2}=\Lambda^2(\R^n)=\{a_{ij}v^i\wedge v^j:v^i,
v^j\in\R^n\}.$$
$\Lambda^2$ is a linear space of dimension
$\frac{n(n-1)}{2}$. If $e_k$ is a normal basis of $\R^n$,
then
$\{e_i\wedge e_j: i<j\}$ is a basis of $\Lambda^2$.
The standard inner product of $\Lambda^2$
is defined by:
$$\lan v_1\wedge v_2,w_1\wedge w_2\ran:=
(v_1\cdot w_1)(v_2\cdot w_2)-(v_1\cdot w_2)(v_2\cdot w_1).$$
So, $\{e_i\wedge e_j\}$ is a normal basis of $\Lambda^2$.

Let $P(\Lambda^2)$ be the projective space getting from
$\Lambda^2$. Recall that there is a nature map $\pi$ from
the unit sphere of $\Lambda^2$ to $P(\Lambda^2)$ which
is a covering map.

Let $\psi$ to be the Pl\"ucker embedding from $G(2,n)$
to $P(\Lambda^2)$, which  endows $G(2,n)$ a Riemannian
metric.
Thus, given a  $b=(e_1,e_2)\in W^{1,2}$, we think of $\varphi(x)=
e_1\wedge e_2$
as a map from $\Omega$ to the unit sphere of $\Lambda^2$ (also
a map to $\Lambda^2$), then the normal of
$\frac{\partial (e_1\wedge e_2)}{\partial x}$ is just the
normal of  $\frac{\partial e_1}{\partial x}\wedge e_2+
e_1\wedge\frac{\partial e_2}{\partial x}$
in $\Lambda^2$. By a direct calculation, we get
$$\begin{array}{lll}
|\frac{\partial (e_1\wedge e_2)}{\partial x}|^2&=&
|\frac{\partial e_1}{\partial x}\wedge e_2+
e_1\wedge\frac{\partial e_2}{\partial x}|^2\\[\mv]
&=&|\frac{\partial e_1}{\partial x}\wedge e_2|^2+
|e_1\wedge\frac{\partial e_2}{\partial x}|^2+2
\lan\frac{\partial e_1}{\partial x}\wedge e_2,
e_1\wedge\frac{\partial e_2}{\partial x}\ran\\[\mv]
&=&|\frac{\partial e_1}{\partial x}|^2+
|\frac{\partial e_2}{\partial x}|^2-2
|e_1\frac{\partial e_2}{\partial x}|^2.
\end{array}$$
So we have
\begin{equation}\label{gaussmap-energy}
|\nabla \varphi|^2=|\nabla b|^2-2|\lan de_1,e_2\ran|^2.
\end{equation}

\vspace{1ex}

Now, we prove (1.1).
Let $(e_1',e_2')$ be a another positively oriented norm basis of $X$. Then we have
$$e_1'=\lambda e_1+\mu e_2,\s e_2'=-\mu e_1+\lambda e_2,$$
where $\lambda=(e_1',e_1)$ and $\mu=(e_1',e_2)$.
We have
$$\frac{\partial e_1'}{\partial x^i}=\frac{\partial\lambda}{\partial x^i}e_1
+\lambda\frac{\partial e_1}{\partial x^i}
+\frac{\partial\mu}{\partial x^i}e_2
+\mu\frac{\partial e_2}{\partial x^i},$$
$$\frac{\partial e_2'}{\partial x^i}=-\frac{\partial\mu}{\partial x^i}e_1
-\mu\frac{\partial e_1}{\partial x^i}
+\frac{\partial\lambda}{\partial x^i}e_2
+\lambda\frac{\partial e_2}{\partial x^i}.$$
We have
$$\begin{array}{lll}
\frac{\partial e_1'}{\partial x^1} \frac{\partial e_2'}{\partial x^2}
&=&-\frac{\partial\lambda}{\partial x^1}
\frac{\partial\mu}{\partial x^2}+\lambda\frac{\partial\lambda}
{\partial x^1}e_1\frac{\partial e_2}{\partial x^2}-
\lambda\mu\frac{\partial e_1}{\partial x^1}
\frac{\partial e_1}{\partial x^2}+\lambda\frac{\partial\lambda}
{\partial x^2}\frac{\partial e_1}{\partial x^1}e_2+\lambda^2
\frac{\partial e_1}{\partial x^1}\frac{\partial e_2}{\partial x^2}\\[\mv]
&&-\mu\frac{\partial \mu}{\partial x^1}e_2\frac{\partial e_1}{\partial x^2}
+\frac{\partial\mu}{\partial x^1}\frac{\partial\lambda}
{\partial x^2}-\mu\frac{\partial\mu}{\partial x^2}\frac{\partial e_2}
{\partial x^1}e_1-\mu^2\frac{\partial e_2}{\partial x^1}
\frac{\partial e_1}{\partial x^2}+\mu\lambda\frac{\partial e_2}
{\partial x^1}\frac{\partial e_2}{\partial x^2},
\end{array}$$
$$\begin{array}{lll}
\frac{\partial e_1'}{\partial x^2} \frac{\partial e_2'}{\partial x^1}
&=&-\frac{\partial\lambda}{\partial x^2}
\frac{\partial\mu}{\partial x^1}+\lambda\frac{\partial\lambda}
{\partial x^2}e_1\frac{\partial e_2}{\partial x^1}-
\lambda\mu\frac{\partial e_1}{\partial x^2}
\frac{\partial e_1}{\partial x^1}+\lambda\frac{\partial\lambda}
{\partial x^1}\frac{\partial e_1}{\partial x^2}e_2+\lambda^2
\frac{\partial e_1}{\partial x^2}\frac{\partial e_2}{\partial x^1}\\[\mv]
&&-\mu\frac{\partial \mu}{\partial x^2}e_2\frac{\partial e_1}{\partial x^1}
+\frac{\partial\mu}{\partial x^2}\frac{\partial\lambda}
{\partial x^1}-\mu\frac{\partial\mu}{\partial x^1}\frac{\partial e_2}
{\partial x^2}e_1-\mu^2\frac{\partial e_2}{\partial x^2}
\frac{\partial e_1}{\partial x^1}+\mu\lambda\frac{\partial e_2}
{\partial x^2}\frac{\partial e_2}{\partial x^1}.
\end{array}$$
We have
$$\begin{array}{lll}
\frac{\partial e_1'}{\partial x^1} \frac{\partial e_2'}{\partial x^2}
-\frac{\partial e_1'}{\partial x^2} \frac{\partial e_2'}{\partial x^1}&=&-2(\frac{\partial\lambda}
{\partial x^1}\frac{\partial\mu}
{\partial x^2}-\frac{\partial\lambda}
{\partial x^2}\frac{\partial\mu}
{\partial x^1})+2(\lambda\frac{\partial\lambda}
{\partial x^1}+\mu\frac{\partial\mu}{\partial x^1})
e_1\frac{\partial e_2}{\partial x^2}
+2(\lambda\frac{\partial\lambda}
{\partial x^2}+\mu\frac{\partial\mu}{\partial x^2})
e_1\frac{\partial e_2}{\partial x^1}\\[\mv]
&&+(\lambda^2+\mu^2)(\frac{\partial e_1}
{\partial x^1}\frac{\partial e_2}{\partial x^2}-
\frac{\partial e_1}{\partial x^2}
\frac{\partial e_2}{\partial x^1}).
\end{array}$$
Since $\lambda^2+\mu^2=1$, we have $\frac{\partial\lambda}
{\partial x^1}\frac{\partial\mu}
{\partial x^2}-\frac{\partial\lambda}
{\partial x^2}\frac{\partial\mu}
{\partial x^1}=0$, and
$\lambda\frac{\partial\lambda}{\partial x^i}+\mu\frac{\partial\mu}
{\partial x^i}=0$, then we get (1.1).\\

We extend $e_1$, $e_2$ to a normal basis $e_3$,
$\cdots$, $e_n\in W^{1,2}$. Such $e_i(i\geq 3)$
exists because $\varphi$ is also a $W^{1,2}$
map from $B$ to $G(n-2,2)$.

We set
$$de_i=w^k_{ij}dx^j\otimes e_k+B^\alpha_{ij} dx^j\otimes
e_\alpha,$$
where $i=1,2$ and $\alpha\in \{3,4,\cdots, n\}$.
Obviously, $w^1_{1i}=w^2_{2i}=0$, $w^1_{2i}=-w^2_{1i}=
\lan \frac{\partial e_1}{\partial x^i},e_2
\ran$, hence (\ref{gaussmap-energy}) is equivalent to
$$|\nabla\varphi|^2=\sum_{ij,\alpha}|B^\alpha_{ij}|^2.$$
 We have
$$\begin{array}{lll}
 \K(\varphi)&=&(w_{11}^ke_k+B_{11}^\alpha n_\alpha)
(w_{22}^ke_k+B_{22}^\alpha n_\alpha)-
(w_{12}^ke_k+B_{12}^\alpha n_\alpha)
(w_{21}^ke_k+B_{21}^\alpha n_\alpha)\\[\mv]
&=&\sum\limits_{\alpha}(B_{11}^\alpha\cdot B_{22}^\alpha-
  |B_{12}^\alpha|^2),
\end{array}$$
therefor we obtain
\begin{eqnarray}
 \K(\varphi)\leq\frac{1}{2}|\nabla\varphi|^2.
\end{eqnarray}
Now, we consider the Gauss map of a conformal map
$f:\Omega\rightarrow\R^n$. Let $v=\frac{1}{2}\log(|\nabla f|^2/2)$
and denote by $X_f$ the Gauss map
induced by $f$.

$X_f$ can be expressed as
$$X_f=(e^{-v}\frac{\partial f }{\partial x^1})\wedge (e^{-v}\frac{\partial f}{\partial x^2}),$$
where $v=\frac{1}{2}\log |\frac{\partial f }{\partial x^1}|^2$. We
will calculate $|\nabla X_f|^2$. Since
$$\frac{\partial^2f }{\partial x^1\partial x^1}\cdot\frac{\partial f }{\partial x^1}
=\frac{1}{2}\frac{\partial }{\partial x^1}\left|\frac{\partial f
}{\partial x^1}\right|^2=e^{2v}\frac{\partial v}{\partial x^1},\s
\frac{\partial^2f }{\partial x^1\partial x^1}\cdot\frac{\partial f }{\partial
x^2}=-\frac{\partial f }{\partial x^1} \cdot \frac{\partial^2
f}{\partial x^1\partial x^2} = -\frac{1}{2}\frac{\partial }{\partial
x^2}\left|\frac{\partial f }{\partial
x^1}\right|^2=-e^{2v}\frac{\partial v }{\partial x^2},$$
and
$$\frac{\partial^2f }{\partial x^1\partial x^1}=A_{11}+\frac{\partial^2f }{\partial x^1\partial x_1}\cdot
 \frac{\partial f }{\partial x^1} e^{-2v}\frac{\partial f }{\partial x^1}
 +\frac{\partial^2f }{\partial x^1}\cdot \frac{\partial f }{\partial x^2}
e^{-2v}\frac{\partial f }{\partial x^2},$$
we get
$$\frac{\partial}{\partial x^1} (e^{-v}\frac{\partial f }{\partial x^1})=
e^{-v}\frac{\partial^2f }{\partial x^1\partial x_1}-e^{-v}\frac{\partial v
}{\partial x^1}\frac{\partial f }{\partial
x^1}=e^{-v}(A_{11}-\frac{\partial f }{\partial x^2}\frac{\partial v
}{\partial x^2}).$$
In the same way, we get
$$\frac{\partial} {\partial x^2}(e^{-v}\frac{\partial f }{\partial
x^1}) =e^{-v}(A_{12}+\frac{\partial f }{\partial x^1}\frac{\partial
v }{\partial x^2}),$$
$$\frac{\partial} {\partial x^1}(e^{-v}\frac{\partial f }{\partial
x^2}) =e^{-v}(A_{21}+\frac{\partial f }{\partial x^2}\frac{\partial
v }{\partial x^1}),$$
$$\frac{\partial} {\partial x^2}(e^{-v}\frac{\partial f }{\partial
x^2})
=e^{-v}(A_{22}-\frac{\partial f }{\partial x^1}\frac{\partial v
}{\partial x^1}).$$
Then, we get
\begin{equation}\label{Xf}
\K(X_f)=e^{-2v}(A_{11}A_{22}-A_{12}^2)=Ke^{2v}
\end{equation}
and
\begin{equation}\label{gauss-A}
|\nabla X_f|^2
=e^{-2v}\sum |A_{ij}|^2,\s i.e.\s |\nabla_{g_f} X_f|^2d\mu_{g_f}
=|A|^2d\mu_{g_f}.
\end{equation}

\subsection{B}
In this part, we will give an alternative proof about that $f_0$ is conformal in theorem 1.12. We need a special case of the following theorem proved by Hardt, Lin and Mou(\cite{HLM}).
\begin{thm}
Let $\Omega$ be a smooth bounded domain in $R^2$, and suppose $1<p<\infty$ and for each $i=1,2,...,$ $u_i\in W^{1,p}(\Omega)$ is a weak solution of
$$div(|\nabla u|^{p-2}\nabla u)+f_i=0$$
with $\sup_i\|u_i\|_{W^{1,p}}+\sup_i\|f_i\|_{L^1}<\infty$. If $u_i\to u$ weakly in $W^{1,p}$, then  $u_i\to u$ strongly in $W^{1,q}$, whenever $1<q<p$.
\end{thm}
Note that $f_m$ is conformal, that is
$$|\frac{\partial f_m}{\partial x}|^2=|\frac{\partial f_m}{\partial y}|^2,\s \frac{\partial f_m}{\partial x}\cdot\frac{\partial f_m}{\partial y}=0. $$
Now by the above theorem we have that
$$f_m\to f_0\s strongly\s in\s W^{1,p},$$
whenever $1<p<2,$
which implies that
$$\frac{\partial f_m}{\partial x}\to \frac{\partial f_0}{\partial x} \s a.e. , \s and\s \frac{\partial f_m}{\partial y}\to \frac{\partial f_0}{\partial y}\s a.e. .$$
Therefor we obtain
$$|\frac{\partial f_0}{\partial x}|^2=|\frac{\partial f_0}{\partial y}|^2,\s \frac{\partial f_0}{\partial x}\cdot\frac{\partial f_0}{\partial y}=0,$$
implying that $f_0$ is conformal.

{}

\vspace{1cm}\sc

YONG LUO

Mathematisches Institut, Albert-Ludwigs-Universit\"at Freiburg,

Eckerstr. 1, 79104 Freiburg, Germany.

{\tt yong.luo@nath.uni-freiburg.de}

\end{document}